\documentclass[11pt, twoside]{amsart}
\usepackage[english]{babel}
\theoremstyle{plain}

\usepackage{amsmath,amsfonts,amsthm,amssymb}
\usepackage{graphicx}
\usepackage{verbatim}
\usepackage{color}
\usepackage{amscd}
\usepackage{verbatim}
\usepackage{mathrsfs}

\addtocounter{page}{1}%
\newtheorem{theorem}{Theorem}[section]

\newtheorem{Thm}[theorem]{Theorem}
\newtheorem{Lem}[theorem]{Lemma}
\newtheorem{Pro}[theorem]{Proposition}

\newtheorem{Rem}[theorem]{Remark}

\newtheorem{Cor}[theorem]{Corollary}

\newtheorem*{corollary*}{Corollary}
\newtheorem*{theorem*}{Theorem}

\def \R {{\mathbb {R}}}

\begin{document}

\setcounter{page}{1}

\newcommand{\M}{{\mathcal M}}
\newcommand{\loc}{{\mathrm{loc}}}
\newcommand{\core}{C_0^{\infty}(\Omega)}
\newcommand{\sob}{W^{1,p}(\Omega)}
\newcommand{\sobloc}{W^{1,p}_{\mathrm{loc}}(\Omega)}
\newcommand{\merhav}{{\mathcal D}^{1,p}}
\newcommand{\be}{\begin{equation}}
\newcommand{\ee}{\end{equation}}
\newcommand{\mysection}[1]{\section{#1}\setcounter{equation}{0}}
\newcommand{\laplace}{\Delta}
\newcommand{\pl}{\laplace_p}
\newcommand{\grad}{\nabla}
\newcommand{\pd}{\partial}
\newcommand{\bo}{\pd}
\newcommand{\csub}{\subset \subset}
\newcommand{\sm}{\setminus}
\newcommand{\ssm}{:}
\newcommand{\diver}{\mathrm{div}\,}
\newcommand{\bea}{\begin{eqnarray}}
\newcommand{\eea}{\end{eqnarray}}
\newcommand{\bean}{\begin{eqnarray*}}
\newcommand{\eean}{\end{eqnarray*}}
\newcommand{\thkl}{\rule[-.5mm]{.3mm}{3mm}}
\newcommand{\cw}{\stackrel{\rightharpoonup}{\rightharpoonup}}
\newcommand{\id}{\operatorname{id}}
\newcommand{\supp}{\operatorname{supp}}
\newcommand{\wlim}{\mbox{ w-lim }}
\newcommand{\mymu}{{x_N^{-p_*}}}
\newcommand{\abs}[1]{\lvert#1\rvert}
\newcommand{\pf}{\noindent \mbox{{\bf Proof}: }}


\renewcommand{\theequation}{\thesection.\arabic{equation}}
\catcode`@=11 \@addtoreset{equation}{section} \catcode`@=12
\newcommand{\Real}{\mathbb{R}}

\newcommand{\B}{\mathbb{B}}

\newcommand{\real}{\mathbb{R}}
\newcommand{\Nat}{\mathbb{N}}
\newcommand{\ZZ}{\mathbb{Z}}
\newcommand{\CC}{\mathbb{C}}
\newcommand{\Pess}{\opname{Pess}}
\newcommand{\Proof}{\mbox{\noindent {\bf Proof} \hspace{2mm}}}
\newcommand{\mbinom}[2]{\left (\!\!{\renewcommand{\arraystretch}{0.5}
\mbox{$\begin{array}[c]{c}  #1\\ #2  \end{array}$}}\!\! \right )}
\newcommand{\brang}[1]{\langle #1 \rangle}
\newcommand{\vstrut}[1]{\rule{0mm}{#1mm}}
\newcommand{\rec}[1]{\frac{1}{#1}}
\newcommand{\set}[1]{\{#1\}}
\newcommand{\dist}[2]{$\mbox{\rm dist}\,(#1,#2)$}
\newcommand{\opname}[1]{\mbox{\rm #1}\,}
\newcommand{\mb}[1]{\;\mbox{ #1 }\;}
\newcommand{\undersym}[2]
 {{\renewcommand{\arraystretch}{0.5}  \mbox{$\begin{array}[t]{c}
 #1\\ #2  \end{array}$}}}
\newlength{\wex}  \newlength{\hex}
\newcommand{\understack}[3]{%
 \settowidth{\wex}{\mbox{$#3$}} \settoheight{\hex}{\mbox{$#1$}}
 \hspace{\wex}  \raisebox{-1.2\hex}{\makebox[-\wex][c]{$#2$}}
 \makebox[\wex][c]{$#1$}   }%
\newcommand{\smit}[1]{\mbox{\small \it #1}}
\newcommand{\lgit}[1]{\mbox{\large \it #1}}
\newcommand{\scts}[1]{\scriptstyle #1}
\newcommand{\scss}[1]{\scriptscriptstyle #1}
\newcommand{\txts}[1]{\textstyle #1}
\newcommand{\dsps}[1]{\displaystyle #1}
\newcommand{\dx}{\,\mathrm{d}x}
\newcommand{\dy}{\,\mathrm{d}y}
\newcommand{\dz}{\,\mathrm{d}z}
\newcommand{\dt}{\,\mathrm{d}t}
\newcommand{\dr}{\,\mathrm{d}r}
\newcommand{\du}{\,\mathrm{d}u}
\newcommand{\dv}{\,\mathrm{d}v}
\newcommand{\dV}{\,\mathrm{d}V}
\newcommand{\ds}{\,\mathrm{d}s}
\newcommand{\dS}{\,\mathrm{d}S}
\newcommand{\dk}{\,\mathrm{d}k}

\newcommand{\dphi}{\,\mathrm{d}\phi}
\newcommand{\dtau}{\,\mathrm{d}\tau}
\newcommand{\dxi}{\,\mathrm{d}\xi}
\newcommand{\deta}{\,\mathrm{d}\eta}
\newcommand{\dsigma}{\,\mathrm{d}\sigma}
\newcommand{\dtheta}{\,\mathrm{d}\theta}
\newcommand{\dnu}{\,\mathrm{d}\nu}

\newcommand{\1}{\mathbf{1}}

\def\ga{\alpha}     \def\gb{\beta}       \def\gg{\gamma}
\def\gc{\chi}       \def\gd{\delta}      \def\ge{\epsilon}
\def\gth{\theta}                         \def\vge{\varepsilon}
\def\gf{\phi}       \def\vgf{\varphi}    \def\gh{\eta}
\def\gi{\iota}      \def\gk{\kappa}      \def\gl{\lambda}
\def\gm{\mu}        \def\gn{\nu}         \def\gp{\pi}
\def\vgp{\varpi}    \def\gr{\rho}        \def\vgr{\varrho}
\def\gs{\sigma}     \def\vgs{\varsigma}  \def\gt{\tau}
\def\gu{\upsilon}   \def\gv{\vartheta}   \def\gw{\omega}
\def\gx{\xi}        \def\gy{\psi}        \def\gz{\zeta}
\def\Gg{\Gamma}     \def\Gd{\Delta}      \def\Gf{\Phi}
\def\Gth{\Theta}
\def\Gl{\Lambda}    \def\Gs{\Sigma}      \def\Gp{\Pi}
\def\Gw{\Omega}     \def\Gx{\Xi}         \def\Gy{\Psi}

\newcommand{\Rn}{\mathbb{R}^n}
 \newcommand{\lanbox}{{\, \vrule height 0.25cm width 0.25cm depth 0.01cm \,}}
  \renewcommand{\qedsymbol}{\lanbox}
   \newcommand{\dd}{\mathrm{d}}
\newcommand{\cqfd}{\begin{flushright}                  
			 $\Box$
                 \end{flushright}}

\renewcommand{\div}{\mathrm{div}}
\newcommand{\red}[1]{{\color{red} #1}}

\pagestyle{headings}
\title{A new eigenvalue problem for free boundary minimal submanifolds in the unit ball}

\author{Baptiste Devyver}
\address{Baptiste Devyver, Department of Mathematics, Technion, 32000 Haifa, Israel}
\email{devyver@technion.ac.il}
\maketitle
\tableofcontents

\section{Introduction}

In this paper, we will study different eigenvalue problems that are related to the index of minimal submanifolds with free boundary in the Euclidean unit ball. Loosely speaking, minimal submanifolds of dimension $k$ with free boundary in the unit ball of $\R^n$ are the critical points for the $k$-volume functional acting on all such submanifolds that lie inside the ball, and whose boundary lies in the unit sphere. It is an interesting problem to understand the interaction between the Morse index (the dimension of the space of perturbation that decrease the area to order two), and the topology of these submanifolds. As an example of this interaction, it was proved recently by P. Sargent \cite{S} and independently by L. Ambrozio, A. Carlotto and B. Sharp \cite{ACS} that the index of a free boundary minimal surface in the unit $3$-ball is at least $\frac{1}{3}(2g+k-1)$, where $g$ is the genus of the surface and $k$ the number of its boundary components. There has been also interest into understanding free boundary minimal surfaces with low index. First, the simplest free boundary minimal surface in the unit ball, namely the flat disk, has index $1$. A. Fraser and R. Schoen \cite{FS1} proved that the index of a non-flat free boundary minimal surface in the unit $3$-ball is at least $3$; later on, it was shown in \cite{D} that actually the index of such a surface is at least $4$. It was shown independently in \cite{D}, \cite{SZ} and \cite{T}, that the index of the so-called {\em critical catenoid}, a minimal annulus with free boundary in the unit $3$-ball introduced by A. Fraser and R. Schoen in \cite{FS1}, is equal to $4$. A natural question is then whether the critical catenoid is the unique (up to congruence) free boundary minimal surface in the unit $3$-ball with index $4$. A motivation for studying this question comes from a celebrated result of F. Urbano \cite{U} asserting the following: firstly, a minimal surface in the $3$-sphere that is not a totally geodesic $S^2$ has index at least $5$, and secondly the Clifford torus is characterized as the unique minimal surface of the $3$-sphere with index equal to $5$.

\medskip

 An easy but crucial technical point in Urbano's proof of the uniqueness of the Clifford torus, is that if the index of a minimal surface in $S^3$ is equal to $5$, then the second eigenvalue of the Jacobi operator has to be equal to $-2$ (with multiplicity $4$). Moreover, for every minimal surface in $S^3$, the $4$ components of the unit normal to the surface are eigenfunctions of the Jacobi operator, associated to the eigenvalue $-2$. Thus, in the case the surface is not an equatorial $S^2$, there is a natural $4$-dimensional space of variations that decrease the area. On a minimal surface with free boundary in the unit $3$-ball, from a spectral point of view, the situation is more complicated: indeed, there are several legitimate eigenvalue problems for the Jacobi operator (e.g. Steklov, Robin, Dirichlet), and each of them is relevant for the study of the index. Indeed, in \cite{D} and \cite{T} the Steklov and Dirichlet spectrum of the critical catenoid are studied in order to compute the index, while in \cite{SZ} the Robin spectrum is used for the same purpose. Nonetheless, while for the critical catenoid the components of the unit normal are Steklov eigenfunctions for the Jacobi operator, this is not the case in general; in \cite{T2}, the critical catenoid is even characterized as the unique non-flat free boundary minimal surface of the unit $3$-ball, for which the components of the unit normal are Steklov eigenfunctions for the Jacobi operator. Furthermore, it does not appear that either the Robin, nor the Steklov eigenvalues and eigenfunctions for the Jacobi operator have any geometric meaning in general. It is thus unclear how to interpret spectrally the information that a free boundary minimal surface has index equal to $4$. Still, it has been proved in \cite{FS2} that the $3$ components of the unit normal to a general free boundary surface in the unit ball of $\R^3$ generate normal variations to the surface that decrease the area. It is however not clear whether these components are eigenfunctions of some eigenvalue problem for the Jacobi operator. 
 
 \medskip

The purpose of the present paper is to show that the components of the unit normal of any minimal surface with free boundary in the unit $n$-ball, are indeed eigenfunctions associated with the eigenvalue $-2$, for some (new) natural eigenvalue problem for the Jacobi operator; this fact has analytic (spectral) consequences for free boundary minimal surfaces in the unit $3$-ball of index $4$. The eigenvalue problem in question can be interpreted as a kind of Steklov eigenvalue problem with an additional non-local term. As previously mentioned, this is in strong analogy with what happens for minimal surfaces in $S^3$. However, let us say right away that unfortunately, we are not able to classify free boundary minimal surfaces in the unit $3$-ball of index $4$, but we hope that our results will be useful in order to make progress towards such a classification. 

\section{Preliminaries}

In all this paper, $\Sigma=\Sigma^k$ denotes a smooth, free boundary $k$-dimensional submanifold in the $n$-dimensional Euclidean ball $\mathbb{B}^n$. As is well-known, such manifolds are precisely the critical points of the $k$-area functional, when deformations are not fixed at the boundary, but rather are only assumed to have boundary included in the unit sphere. One can consider the quadratic form associated with the second variation of $k$-area of $\Sigma$. Restricted to normal deformations, it gives rise to a quadratic form $Q$ on the normal bundle $\mathscr{N}$ over $\Sigma$, defined as follows:

$$Q(W)=\int_\Sigma(||D^\perp W||_{HS}^2-|(\sigma(\cdot,\cdot),W)|^2)-\int_{\partial\Sigma}||W||^2.$$
Here, $D^\perp$ is the connection on the normal bundle, obtained by orthogonal projection of the connection of $\R^n$ onto $\mathscr{N}$, and $||D^\perp W||_{HS}^2$ denotes the Hilbert-Schmidt norm: for $(e_i)_{i=1,\cdots,k}$ an orthonormal basis of $T_p\Sigma$, by definition

$$||D^\perp W||_{HS}^2(p)=\sum_{i=1}^k||D^\perp_{e_i} W||^2.$$
Also, $\sigma$ denotes the (vector-valued) second fundamental form $\sigma(X,Y)=(\nabla_XY)^\perp$, and

$$|(\sigma(\cdot,\cdot),W)|^2=\sum_{i,j=1}^k|(\sigma(e_i,e_j),W)|^2.$$
For $W:\Sigma\to \R^n$, we will denote by $W^\perp$ the projection of $W$ onto $\mathscr{N}$, and by $W^\top$ the projection of $W$ onto the tangent bundle $T\Sigma$. The stability operator $J$ is the second-order differential operator naturally associated to $Q$; it satisfies

$$Q(W)=\int_\Sigma (JW,W),$$
for all $W\in \Gamma(\mathscr{N})$ vanishing at the boundary of $\Sigma$. Explicitly,

$$JW=\Delta_\Sigma^\perp W-\mathscr{S}(W),$$
with $\Delta_\Sigma^\perp$ the Laplacian on the normal bundle, defined locally in an orthonormal basis  $\{e_i\}_{i=1,2}$ of $T_p\Sigma$ by

$$\Delta_\Sigma^\perp W=-\sum_{i=1}^k(\nabla_{e_i}\nabla_{e_i}W)^\perp+\sum_{i=1}^k (\nabla_{(\nabla_{e_i}e_i)^T}W)^\perp.$$
Moreover, $\mathscr{S}(W)$, the Simons operator, is defined locally in an orthonormal basis  $(e_i)_{i=1,\cdots,k}$ of $T_p\Sigma$ by

$$\mathscr{S}(W)=\sum_{i=1}^k (\sigma(e_i,e_j),W)\sigma(e_i,e_j).$$
The {\em index} of $\Sigma$ as a free boundary minimal surface of $\mathbb{B}^n$, which will be denoted $\mathrm{ind}(\Sigma)$, is by definition the maximal dimension of a subspace $V$ of $C^\infty(\mathscr{N})$, on which $Q$ is negative. 

For $v\in\R^n$, one denotes $v^\perp$ the section of the normal bundle $\mathcal{N}$ obtained from $v$ by orthogonal projection. In \cite{FS2}, it was shown that for every $v\in \mathbb{R}^n$

\begin{equation}\label{FSN}
Q(v^\perp,v^\perp)=-k\int_\Sigma |v^\perp|^2,
\end{equation}
which implies that $\mathrm{ind}(\Sigma)\geq n$, unless $\Sigma$ is contained in a subspace of dimension $n-1$. In \cite{D}, this lower bound was improved in the case $k=2$, $n=3$, and it was shown there that for $\Sigma^2\subset \mathbb{B}^3$ orientable which is not a flat disk, $\mathrm{ind}(\Sigma)\geq 4=3+1$. Furthermore, this inequality is sharp, since it was proven there that the so-called critical catenoid has index precisely equal to $4$ (see also \cite{SZ}, \cite{T}).

\section{A slight variation on Urbano's index characterization of the Clifford torus}

In this section, we present a variation on the proof, due to F. Urbano \cite{U}, of the celebrated index characterization of the Clifford torus as the unique closed minimal surface in $S^3$ having index $5$. This section is independent of the rest of the paper, and serves as a motivation for the results that will be obtained in the next sections. So, we let $\Sigma$ be a minimal surface of index $5$ in $S^3$ of genus $g$, and we assume that $\Sigma$ is not a totally geodesic sphere $S^2$. One first wants to estimate the genus $g$. Let $J=\Delta-|A|^2-2$ be the Jacobi operator. It is easily checked that $JN=-2N$, where $N$ is the unit normal to $\Sigma$ in $S^3$. As in \cite{U}, since $\Sigma$ is not an equatorial sphere, the vector space $Span\{N_i\,;\,i=1,\cdots,4\}$ is $4$-dimensional, hence the first eigenvalue of $J$, $\lambda_1(J)$, is strictly less than $-2$. Let $\rho>0$ be an associated first eigenfunction. Since the index is exactly $5$, the second eigenvalue of $J$ must be $-2$, hence by the min-max characterization of the eigenvalues one gets

\begin{equation}\label{-2_eigen}
S(u,u)\geq -2\int_\Sigma u^2,
\end{equation}
for all $u\in C^\infty(\Sigma)$ such that $\int_\Sigma u\rho=0$. Here, the quadratic form $S$ is defined by

$$S(u,u)=\int_\Sigma |\nabla u|^2-|A|^2u^2-2u^2.$$
Let $\varphi:\Sigma\to S^2$ be a conformal map of degree $\leq \left[\frac{g+3}{2}\right]$ (such a map exists by general results from algebraic geometry, see \cite[Theorem 4]{R}). By a well-known lemma of P. Li and S. T. Yau (see the proof of \cite[Theorem 1]{LY}), up to composing $\varphi$ by a Mobius transformation, one can assume that $\varphi$ is balanced, that is 

$$\int_\Sigma \varphi\rho=0.$$
We plug the coordinates $\varphi_i$ of $\varphi$ into \eqref{-2_eigen}, and sum over $i=1,\cdots,4$. Then, one obtains

$$\int_\Sigma|\nabla \varphi|^2\geq\int_\Sigma |A|^2.$$
By the Gauss equations, $|A|^2=2-2K$, and using the Gauss-Bonnet theorem, one gets

$$\int_\Sigma |A|^2=2A(\Sigma)+8\pi(g-1).$$
On the other hand, by the conformality of $\varphi$,

$$\int_\Sigma |\nabla \varphi|^2=2\,deg(\varphi)\,A(\Sigma)\leq 8\pi \left[\frac{g+3}{2}\right].$$

\medskip
CLAIM 1: the area of $\Sigma$ is strictly greater than $4\pi$.

\medskip
Assuming the result of the claim for the moment, one finds that $g< \left[\frac{g+3}{2}\right]$, and from this one concludes easily that $g$ is equal to zero or one. By a result of Almgren \cite{A}, if $g=0$ then $\Sigma$ is a totally geodesic sphere, which is excluded, so $\Sigma$ is topologically a torus. At this point one could resort to the solution of the Lawson conjecture by S. Brendle \cite{B}, to conclude that $\Sigma$ must be the Clifford torus. However, there is an alternative, more elementary argument, which relies on the following:

\medskip
CLAIM 2: the first non-zero eigenvalue of the Laplacian on $\Sigma$ is equal to $2$.

\medskip
This immediately implies that $\Sigma$ is the Clifford torus, according to \cite[Theorem 4]{MR}.

\medskip

{\em Proof of Claim 1}

\medskip

According to \cite{LY}, Proposition 1 and Fact 2, $A(\Sigma)\geq 4\pi$. Hence, it is enough to exclude the equality case. This follows from arguments in \cite{MR}. According to (1.12) therein, in the equality case, there exists a constant unit vector $g\in S^3$ such that $g^N$ (the normal part of $g$ along $\Sigma$) vanishes identically. The discussion in \cite{MR}, proof of Theorem 1, then implies that $\Sigma$ must be the standard $2$-sphere. This is excluded by Almgren's result.

\medskip

{\em Proof of Claim 2}

\medskip

The proof follows the proof of \cite[Prop. 6.2]{D}. Assume by contradiction that there is a non-constant function $h$ on $\Sigma$ with $\Delta h=\lambda h$ for some $\lambda<2$. Since $\Sigma$ is minimal, $\Delta_\Sigma x_i=2x_i$ for every $i=1,2,3,4$. Let $\mathcal{V}$ be the span of $\{x_i\}_{i=1}^4$, of the constant function $1$, and of $h$. Then, $\mathcal{V}$ is $6$-dimensional, otherwise $\Sigma$ is (contained in) a $2$-sphere inside $S^3$, which is excluded. Denote by $R(u,u)$ the quadratic form

$$R(u,u)=\int_\Sigma |\nabla u|^2-2u^2,$$
then $R\leq 0$ on $\mathcal{V}$. Observe that since $g\neq 0$, $|A|$ has only isolated zeroes: indeed, this follows from Gauss' equation, which implies that the zeroes of $|A|$ are precisely the points where the Gauss curvature of $\Sigma$ is equal to one, and the fact that these are isolated (see \cite[Lemma 1.4]{L}). Hence, it follows that for any smooth function $u$ on $\Sigma$,

$$S(u,u)<R(u,u).$$
Hence, $S$ is negative definite in restriction to $\mathcal{V}$, and the index is at least $6$, a contradiction.

\cqfd

\bigskip

Quite naturally, one is led to ask whether the same approach works, in order to characterize unique free boundary, orientable minimal surfaces of $\mathbb{B}^3$ of index $4$. Conjecturaly, the critical catenoid is the only such surface, up to congruence. Notice that most steps in the above proof of Urbano's theorem easily adapt to the case of free boundary minimal surfaces: first, it is known that if $\Sigma$ is a free boundary minimal surface in $\mathbb{B}^3$ which is not a flat disk, then $2A(\Sigma)=L(\partial\Sigma)> 2\pi$ (see \cite{K}), and if in addition $\Sigma$ is topologically an annulus such that the first non-zero Steklov eigenvalue for the Laplacian is equal to $1$, then $\Sigma$ is congruent to the critical catenoid (see \cite{FS2}). It is also known that if $\Sigma$ has index $4$, then the first non-zero Steklov eigenvalue for the Laplacian is equal to $1$ (see \cite{D}). In fact, the only ingredient that is crucially missing is an inequality such as \eqref{-2_eigen}, as well as good test functions to be fed in it.

\section{An eigenvalue problem on the critical catenoid}

As explained in Section 3, a key ingredient in Urbano's proof is the following inequality, for $\Sigma\hookrightarrow S^3$ minimal with index 5: for $u\in C^\infty(\Sigma),$

\begin{equation}\label{ineq_Ur}
Q(u,u)\geq -2\int_\Sigma u^2,
\end{equation}
provided $\int_\Sigma \rho u=0$, where $\rho$ is a first eigenfunction for the Jacobi operator. The inequality \eqref{ineq_Ur} is equivalent to the fact that $-2$ is the second eigenvalue of the Jacobi operator, which is necessarily the case if the index is $5$. The eigenfunctions of $J$ associated to the eigenvalue $-2$ are precisely the four components of the normal vector $N$ to $\Sigma$. In the case of free boundary minimal surfaces in $\mathbb{B}^3$, according to \eqref{FSN} it is natural to expect that the components of the normal vector will play a role, and that the value $-2$ should appear. By analogy, one could expect an inequality such as \eqref{ineq_Ur} to hold true. Unfortunately, as we shall demonstrate now, the obvious generalization of \eqref{ineq_Ur} does not hold on the critical catenoid itself.

As is well-known (see \cite[Section 5]{CFP}), the index of a free boundary minimal surface in $\B^3$ is obtained as the number of negative eigenvalues of the following Robin boundary problem for the Jacobi operator:

\begin{equation}\label{BVPg}
\left\{
\begin{array}{lcr}
Ju=\gamma u\mbox{ on }\Sigma\\
\frac{\partial u}{\partial\nu}= u\mbox{ on }\partial\Sigma.
\end{array}
\right.
\end{equation}
Let us denote $\gamma_1\leq \gamma_2\leq\cdots\leq\gamma_n\leq\cdots$ the associated spectrum. There is also a variational (min-max) characterization of the eigenvalues $\gamma_i$. For example, one has

$$\gamma_1=\inf_{u\neq 0}\frac{Q(u)}{\int_\Sigma u^2}.$$
From this characterization of $\gamma_1$ and the Harnack inequality, one concludes by standard arguments that $\gamma_1$ is simple and that the first eigenfunction $\rho_1$ is positive in the interior of $\Sigma$. One is lead to ask the question: assuming that the index of $\Sigma$ is equal to $4$, does the inequality

\begin{equation}\label{ineq_Ur2}
Q(u)\geq -2\int_\Sigma u^2,
\end{equation}
hold, provided $\int_\Sigma \rho_1 u=0$? Equivalently, is it true that $\gamma_2\geq 2$? We are going to show by elementary arguments that the answer to this question is negative, even in the case of the critical catenoid itself, and even more:

\begin{Pro}\label{energy-cat}

 On the critical catenoid $\Sigma\subset \mathbb{B}^3$, the boundary value problem \eqref{BVPg} admits $4$ negative eigenvalues that are strictly less  than $2$.

\end{Pro}
This implies that for the critical catenoid, in order to ensure that the inequality $Q(u,u)\geq -2\int_\Sigma u^2$ holds, one needs to impose {\em four}, and not just one, orthogonality conditions on $u$. In fact, since the index of the critical catenoid is $4$, one already knows that $Q(u,u)\geq0$ if $u$ is orthogonal to the first four eigenfunctions of \eqref{BVPg}, and Proposition \ref{energy-cat} tells us that one cannot do better, even if one is interested in the weaker inequality \eqref{ineq_Ur}. However, $-2$ is the crucial value in an inequality such as \eqref{ineq_Ur}, that allows one to control the topology of $\Sigma$ by applying it to (the components of) a conformal map $\varphi$ from $\Sigma$ into $S^2$: indeed, this comes from the fact that according to \cite[Theorem 5.4]{FS1},

$$\sum_{i=1}^3\int_{\partial\Sigma}\varphi_i^2=L(\partial \Sigma)=2A(\Sigma)=2\sum_{i=1}^3 \int_{\Sigma} \varphi^2.$$
Hence, Urbano's proof cannot be extended straightforwardly to characterize free boundary minimal surfaces of $\B^3$ of index $4$. In the next two sections, we shall present an alternative inequality, weaker than \eqref{ineq_Ur2}, that holds for any free boundary minimal surface in $\mathbb{B}^3$ with index $4$, and which we hope should play the role of\eqref{ineq_Ur2}, in order to characterize such surfaces.\\

\noindent{\em Proof of Proposition \ref{energy-cat}:} we claim that the first eigenfunction $\rho_1$ of \eqref{BVPg} has the symmetries of the critical catenoid. We recall (see \cite{D}) that the critical catenoid can be parametrized by

$$X(s,\theta)=a(\cosh(s)\cos(\theta),\cosh(s)\sin(\theta),s),\,s\in [-T,T],\,\theta\in [0,2\pi],$$
where $T$ is the unique positive solution of $T\tanh(T)=1$, and $a=(T\cosh(T))^{-1}$. We claim that $\rho_1=\rho_1(s)$ and is an even, positive function. To prove this, let us denote by $\Theta$ the group of symmetries of the catenoid, i.e. the group generated by the reflection w.r.t. the $xy$-plane, and rotations around the $z$-axis. Then, for every $g\in\Theta$, $\rho_1\circ g$ solves the same boundary value problem as $\rho_1$. Averaging the functions $\rho_1\circ g$ over $\Theta$ with its Haar measure, we obtain a positive, $\Theta$-invariant eigenfunction associated to $\gamma_1$. Since $\gamma_1$ is simple, this function must be a constant positive multiple of $\rho_1$, thus $\rho_1$ is $\Theta$-invariant, and the claim follows. 

Notice now that if $u$ is solution of \eqref{BVPg}, then for every $w$,

$$Q(u,w)=\gamma\int_\Sigma u w.$$
Thus, since $\rho_1$ is even and $v_z^\perp$ is odd (as functions of $s$),

\begin{equation}\label{orthog1}
Q(\rho_1,v_z^\perp)=\gamma_1\int_\Sigma \rho_1v_z^\perp=0.
\end{equation}
Actually, since $\rho_1$ is ``radial'', it is also first eigenfunction for a boundary value problem involving the ``radial part'' of the Jacobi operator $\mathscr{L}_0$. To be more precise, let us recall (see \cite{D}) the expression of the Jacobi operator $J$ in $(s,\theta)$ coordinates: 

$$Ju(s,\theta)=-\frac{1}{a^2\cosh^2(s)}\left(\frac{\partial^2 u}{\partial s^2}+\frac{\partial^2u}{\partial \theta^2}\right)-\frac{2}{a^2\cosh^4(s)}u.$$
Functions on $\Sigma$ depending only on the variable $s$ will be called radial. Let us define a ``radial'' Jacobi operator

$$\mathscr{L}_0:=-\frac{1}{a^2\cosh^2(s)}\frac{\partial^2}{\partial s^2}-\frac{2}{a^2\cosh^4(s)},$$
and let us consider the following Robin eigenvalue problem for $\mathscr{L}_0$:
\begin{equation}\label{BVPrad}
\left\{
\begin{array}{lcr}
\mathscr{L}_0\,u(s)=\gamma u(s),\,\,s\in(-T,T)\\
\frac{d u}{ds}=\pm \frac{1}{T} u,\,\,s=\pm T.
\end{array}
\right.
\end{equation}
The eigenvalues of \eqref{BVPrad} have a variational (min-max) characterization in terms of the Rayleigh quotients

$$\frac{Q(u(s))}{\int_\Sigma u(s)^2}$$
for radial functions. In particular, the second eigenvalue of \eqref{BVPrad} is given by

$$\inf \frac{Q(u(s))}{\int_\Sigma u^2},$$
where the infimum is taken over all radial functions $u\not\equiv0$ satisfying $\int_\Sigma u\rho_1=0$. The function $v_z^\perp=\tanh(s)$ is radial, and given that it is odd and $\rho_1$ is even, it follows that

$$\int_\Sigma v_z^\perp\rho_1=0.$$
Therefore, the second eigenvalue of \eqref{BVPrad} is less or equal to

$$\frac{Q(v_z^\perp)}{\int_\Sigma |v_z^\perp|^2}=-2,$$
with equality if and only if $v_z^\perp$ is (second) eigenfunction for \eqref{BVPrad}. Since $Jv_z^\perp=0\neq -2 v_z^\perp$, this cannot be true, and therefore the second eigenvalue of \eqref{BVPrad} has to be $<-2$. Consequently, \eqref{BVPrad} has two eigenvalues $<-2$.\\

Now, let us define a quadratic form $Q_1$ on radial functions by:

$$Q_1(a(s))=Q(a(s)\cos(\theta))=Q(a(s)\sin(\theta)).$$
It is naturally associated to the radial operator 

$$\mathscr{L}_1=\mathscr{L}_0+\frac{1}{a^2\cosh^2(s)}.$$
Let us consider the following Robin boundary value problem for $\mathscr{L}_1$:

\begin{equation}\label{BVPcos}
\left\{
\begin{array}{lcr}
\mathscr{L}_1\,u(s)=\gamma u(s),\,\,s\in(-T,T)\\
\frac{d u}{ds}=\pm \frac{1}{T} u,\,\,s=\pm T.
\end{array}
\right.
\end{equation}
Recall that 

$$v_x^\perp=\Lambda(s)\cos(\theta),\,\,v_y^\perp=\Lambda(s)\sin(\theta),$$
for some explicit (even) radial function $\Lambda(s)$. By the characterization of the first eigenvalue of \eqref{BVPcos} in terms of Rayleigh quotients, the first eigenvalue of \eqref{BVPcos} is less or equal to

$$\frac{Q_1(\Lambda)}{\int_\Sigma \Lambda(s)^2\cos^2(\theta)}=\frac{Q(v_x^\perp)}{\int_\Sigma |v_x^\perp|^2}=\frac{Q(v_y^\perp)}{\int_\Sigma |v_y^\perp|^2}=\frac{Q_1(\Lambda)}{\int_\Sigma \Lambda(s)^2\sin^2(\theta)}=-2,$$
and equality holds if and only if $\Lambda(s)$ is the first eigenfunction of \eqref{BVPcos}. However, since $\mathscr{L}_1(\Lambda)=Jv_x^\perp=0\neq -2\Lambda(s)$, equality cannot hold. So, there is a first eigenfunction $\tilde{\rho}_1$ associated to a first eigenvalue $\tilde{\gamma}_1<-2$ of \eqref{BVPcos}. Thus, we get two eigenfunctions $\tilde{\rho}_1\cos(\theta)$, $\tilde{\rho_1}\sin(\theta)$ of \eqref{BVPg}, associated to the eigenvalue $\tilde{\gamma}_1<-2$. So, we have obtained $4$ eigenvalues of \eqref{BVPg}, that are $<-2$, and such that the associated eigenfunctions are linearly independent. The result follows.

\cqfd








\section{A spectral problem}

In all this section, $\Sigma=\Sigma^k$ denotes a free boundary minimal submanifold of dimension $k$ in $\mathbb{B}^n$. The normal vector bundle $\mathscr{N}$ carries a natural metric induced by the metric of $\mathbb{R}^n$, and a natural inner product:

$$(W_1,W_2)_{L^2}:=\int_\Sigma (W_1(p),W_2(p))\,dvol(p).$$
We will denote by $L^2(\mathscr{N})$ the Hilbert space of $L^2$ sections of $\mathscr{N}$. One defines the Sobolev space $W^{1,2}(\mathscr{N})$, consisting of all sections $W\in L^2(\mathscr{N})$, such that 

$$\int_\Sigma ||(D^\perp W)(p)||_{HS}^2\,dvol(p)<\infty.$$
Let us consider the following two spaces of sections of the normal bundle:

$$\mathscr{H}=\{W\in L^2(\mathscr{N})\,;\,JW=0\hbox{ in }int(\Sigma)\},$$
and
$$\mathscr{E}=\{W\in W^{1,2}(\mathscr{N})\,;\,JW=0\hbox{ in }int(\Sigma)\}\subset \mathscr{H},$$
where the equation $JW=0$ is intended in the weak sense. By elliptic regularity, 

$$\mathscr{H}\subset C^\infty(int(\Sigma)).$$
In fact, elliptic regularity implies that $\mathscr{H}$ is a \textit{closed} subspace of $L^2(\mathscr{N})$, hence a Hilbert space if it is endowed with the $L^2$ inner product. By analogy with complex analysis, $\mathscr{H}$ or $\mathscr{E}$ could be called {\em Hardy spaces} of normal sections. Analogously, $\mathscr{E}$ is a closed subspace of $W^{1,2}(\mathscr{N})$, hence a Hilbert space for the $||\cdot||_{W^{1,2}}$ norm. As follows from the trace theorem,

$$\mathscr{E}\subset W^{1/2,2}(\partial \Sigma)\subset L^2(\partial\Sigma).$$
Hence, the quadratic form $Q$ is naturally defined on $\mathscr{E}$, endowed with the $W^{1,2}(\mathscr{N})$-norm. 

\begin{Pro}\label{SA}

There exists a self-adjoint operator $\mathscr{A}$ with $\mathscr{E}=\mathcal{D}(\mathscr{A}^{1/2})$, such that, for every $W_1$ and $W_2$ in $\mathscr{E}$,

$$Q(W_1,W_2)=\int_\Sigma (W_1,\mathscr{A}W_2)=\int_\Sigma(\mathscr{A}W_1,W_2).$$
Furthermore,  the spectrum of $\mathscr{A}$ consists in a discrete sequence $\mu_0\leq \mu_1\leq\cdots\leq \mu_k\leq\cdots$, tending to $+\infty$.

\end{Pro}

\begin{Rem}
{\em 
The operator $\mathscr{A}$ given by Proposition \eqref{SA} is a {\em non-local} operator.
}
\end{Rem}

\begin{proof} By spectral theory (see \cite[Theorem VIII.15]{RS}), the existence of $\mathscr{A}$ self-adjoint associated with the quadratic form $Q$ follows from the quadratic form $Q$ being bounded from below and closed, in restriction to $\mathscr{E}$. Since $\mathscr{E}$ is closed in $W^{1,2}(\mathscr{N})$, this is equivalent to the existence of $\omega\geq0$, $\alpha>0$, such that, for every $W\in\mathscr{E}$,

\begin{equation}\label{elliptic1}
Q(W)+\omega ||W||_2^2\geq \alpha ||W||_{W^{1,2}}^2.
\end{equation}
We will prove below that \eqref{elliptic1} holds, hence the existence of $\mathscr{A}$ self-adjoint; but for now, let us explain why the statement on the eigenvalues of $\mathscr{A}$ holds. Clearly, it follows from the fact that the resolvent of $\mathscr{A}$ is compact. It is a well-known general fact that provided \eqref{elliptic1} holds, the resolvent of $\mathscr{A}$ is compact if and only if the embedding $\mathscr{E} \hookrightarrow L^2(\mathscr{N})$ is compact. For the sake of completeness, let us detail this point. First, inequality \eqref{elliptic1} implies that $(\mathscr{A}+\omega Id)^{-1/2}$ is bounded from $L^2$ to $\mathscr{E} \subset L^2$. Then, given that $(\mathscr{A}+\omega Id)^{-1/2}$ is bounded on $L^2$, by using the Spectral Theorem for self-adjoint, compact operators one concludes easily that the resolvent $(\mathscr{A}+\omega Id)^{-1}$ being compact is equivalent to $(\mathscr{A}+\omega Id)^{-1/2}$ being compact. Note now that the converse inequality to \eqref{elliptic1},

$$Q(W)+\omega ||W||_2^2\leq \beta ||W||_{W^{1,2}}^2,$$
trivially holds with $\beta=\omega+1$, and implies that $(\mathscr{A}+\omega Id)^{1/2}$ is bounded from $\mathscr{E}$ to $L^2$. Writing the inclusion $\iota : \mathscr{E}\hookrightarrow L^2$ as

$$\iota = (\mathscr{A}+\omega Id)^{-1/2}(\mathscr{A}+\omega Id)^{1/2}:\mathscr{E}\to L^2\to L^2,$$
one sees that $(\mathscr{A}+\omega Id)^{-1/2}$ being compact implies that $\iota$ is compact. In the other direction, writing

$$(\mathscr{A}+\omega Id)^{-1/2}: L^2\to \mathscr{E}\hookrightarrow L^2,$$
where the first arrow is bounded and the last arrow is $\iota$, one sees that $\iota$ compact implies $(\mathscr{A}+\omega Id)^{-1/2}$ compact. Thus, the point is proved.

In our particular situation, by Rellich's theorem, $W^{1,2}(\mathscr{N})\hookrightarrow L^2(\mathscr{N})$ is compact, and since $\mathscr{E}\hookrightarrow W^{1,2}(\mathscr{N})$ is bounded by definition of the norm on $\mathscr{E}$, one concludes that, provided that \eqref{elliptic1} holds, the spectrum of $\mathscr{A}$ consists in a discrete sequence $\mu_0\leq \mu_1\leq\cdots\leq \mu_k\leq\cdots$, tending to $+\infty$.

It remains to prove \eqref{elliptic1}. In fact we will prove the following stronger inequality: there exists $\omega\geq0$, $\alpha>0$, such that, for every $W\in W^{1,2}(\mathscr{N})$,

\begin{equation}\label{elliptic2}
Q(W)+\omega ||W||_2^2\geq \alpha ||W||_{W^{1,2}}^2.
\end{equation}
The inequality \eqref{elliptic2} with $\alpha=\frac{1}{2}$ is an easy consequence of the fact that the Simons operator is bounded, and of the following inequality: there exists $C>0$ such that, for every $W\in W^{1,2}(\mathscr{N})$,

\begin{equation}\label{elliptic3}
\int_{\partial\Sigma} ||W||^2\leq \frac{1}{2}\int_\Sigma ||D^\perp W||^2_{HS}+C\int_\Sigma ||W||^2,
\end{equation}
which is a consequence of \cite[Lemma 2.3]{AM} with the choices $T:W^{1,2}\to L^2$ the ``restriction to the boundary'' operator, $S=id:W^{1,2}(\mathscr{N})\to L^2(\mathscr{N}|_{\partial \Sigma})$ and $\varepsilon=\frac{1}{2}$ therein.

\end{proof}

\section{The $v^\perp$ are eigenfunctions, and consequences}

In this section, we keep the notations that have been defined in Section 5.

\begin{Thm}\label{pseudo}

Let $\Sigma^k\subset \B^n$ be an immersed, free boundary minimal submanifold of dimension $2\leq k\leq n-1$. For every $v\in \R^n\setminus\{0\}$, the normal vector field $v^\perp$ on $\Sigma^k$, defined as the orthogonal projection of $v$ on the normal bundle of $\Sigma^k$, is an eigenfunction of $\mathscr{A}$, associated to the eigenvalue $-k$. In other words, for every normal vector field $W\in \mathscr{E}\subset{W}^{1,2}(\mathscr{N})$, i.e. $JW=0$ in $int(\Sigma)$, then

$$Q(v^\perp,W)=-k\int_\Sigma(v^\perp,W).$$
\end{Thm}
The proof of Theorem \ref{pseudo} follows from the following formula, which we think might be of independent interest:

\begin{Pro}\label{pro:JW}
l
Let $\Sigma^k\subset \B^n$ be an immersed, free boundary minimal submanifold of dimension $2\leq k\leq n-1$. For every $v\in \R^n\setminus\{0\}$ and every normal vector field $W\in \Gamma(\mathscr{N})$, there holds:

\begin{equation}\label{IPP}
Q(v^\perp,W)=-k\int_\Sigma (v^\perp, W)+\int_\Sigma (JW,(v,x)x+\frac{1}{2}(1-|x|^2)v).
\end{equation}

\end{Pro}
Note that the formula \eqref{IPP} implies in particular that

$$Q(v^\perp)=-k\int_\Sigma |v^\perp|^2,$$
a formula that is already known from the work of A. Fraser and R. Schoen (see \cite[Theorem 3.1]{FS2}). However, our proof is different from theirs. The new ingredient, and the main trick in the proof, is to use the properties of the carefully chosen vector field $(v,x)x+\frac{1}{2}(1-|x|^2)v$.

\begin{proof}

Let $v\in\mathbb{R}^n\setminus \{0\}$, and let $Y$ be the (non-tangential) vector field on $\Sigma$ defined by

$$Y=(v,x)x+\frac{1}{2}(1-|x|^2)v.$$
We will denote by $Y^\perp$ the orthogonal projection of $Y$ on the normal bundle of $\Sigma$. By integration by parts, we have

$$\int_\Sigma (JW,Y)=\int_\Sigma (JW, Y^\perp)=\int_\Sigma (W,JY^\perp)+\int_{\partial \Sigma} \left(W,\frac{\pd Y^\perp}{\pd \nu}\right)-\left(\frac{\partial W}{\partial \nu},Y^\perp\right).$$
Note that $Y^\perp \equiv0$ at the boundary of $\Sigma$, therefore, remembering that $W$ is normal at every point, one obtains

\begin{equation}\label{eq:1}
\int_\Sigma (JW,Y)=\int_\Sigma (W,JY^\perp)+\int_{\partial \Sigma} \left(W,\nabla_{\nu}^\perp Y^\perp\right).
\end{equation}
On the other hand, since $Jv^\perp=0$, one has

$$Q(v^\perp,W)=\int_{\partial \Sigma}\left(W,\frac{\pd v^\perp}{\pd \nu}-v^\perp\right),$$
and since $W$ is orthogonal to $\Sigma$,

\begin{equation}\label{eq:2}
Q(v^\perp,W)=\int_{\partial \Sigma}\left(W,D_\nu^\perp v^\perp-v^\perp\right).
\end{equation}
We claim that the following two identities hold:

\begin{equation}\label{eq:3}
JY^\perp=kv^\perp,
\end{equation}
and

\begin{equation}\label{eq:4}
\nabla_{\nu}^\perp Y^\perp=D_\nu^\perp v^\perp-v^\perp
\end{equation}
Clearly, \eqref{eq:3} and \eqref{eq:4}, together with \eqref{eq:1} and \eqref{eq:2} imply the result of Proposition \ref{pro:JW}. In order to prove these, we use the following easy computational lemma, whose proof is postponed for the moment:

\begin{Lem}\label{lem:tech}

Let $p$ be a point of $\Sigma$. Let $(e_1,\cdots,e_k)$ be a local orthonormal basis of $T\Sigma$ around $p$, and let $(N_{k+1},\cdots,N_n)$ be a local orthonormal basis of $\Gamma(\mathscr{N})$ around $p$. Let us denote by $D_i$ the covariant derivative $D_{e_i}$ of $\R^n$ in the direction $e_i$, and by $D^\perp_i$ its the projection onto the normal bundle of $\Sigma$. We assume that at the point $p$, $D_{e_i}e_j$ is normal to $\Sigma$, and that $D^\perp N_i=0$ for all $i=k+1,\cdots,n$. Denote also by $x^\top$ and $v^\top$ the orthogonal projection of $x$ and $v$ respectively onto the tangent space to $\Sigma$. Then, one has at the point $p$:

\begin{itemize}

\item[(i)] $D_i^\perp x^\perp =-\sigma(x^\top,e_i)$.

\item[(ii)] $D_i^\perp v^\perp=-\sigma(v^\top,e_i)$.

\item[(iii)] $\Delta_\Sigma |x|^2=-2k$ (where $\Delta_\Sigma$ is the Laplacian on $\Sigma$).

\item[(iv)] $D_\nu^\perp x^\perp=-\sigma(\nu,\nu)$.

\item[(v)] $D_\nu^\perp x^\perp=-(v,\nu)\sigma(v^\top,\nu)$.

\end{itemize}

\end{Lem}
With this lemma at hand, let us conclude the proof. Using that $Jv^\perp=0$ and (iii) of Lemma \ref{lem:tech}, one has

$$\begin{array}{rcl}
JY^\perp&=&J\left((v,x)x^\perp+\frac{1}{2}(1-|x|^2)v^\perp\right)\\\\
&=& \left(\Delta_\Sigma (v,x)\right)x^\perp+(v,x)Jx^\perp-2\sum_{i=1}^k(v,D_ix)D_i^\perp x^\perp+kv^\perp +\sum_{i=1}^k (D_i|x|^2)D_i^\perp v^\perp.
\end{array}$$
As is well-known, the normal vector field $x^\perp$ is a Jacobi vector field, that is:

$$Jx^\perp=0,$$
and furthermore, since $x$ is harmonic (by minimality of $\Sigma$), one obtains by using (i) and (ii) of Lemma \ref{lem:tech} that

$$\begin{array}{rcl}
JY^\perp&=&kv^\perp -2\sum_{i=1}^k(v,D_ix)D_i^\perp x^\perp +\sum_{i=1}^k (D_i|x|^2)D_i^\perp v^\perp\\\\
&=& kv^\perp -2\sum_{i=1}^k v_iD_i^\perp x^\perp +2\sum_{i=1}^k x_iD_i^\perp v^\perp\\\\
&=& kv^\perp +2\sigma(x^\top,v^\top)-2\sigma(x^\top,v^\top)\\\\
&=& kv^\perp.
\end{array}$$
This concludes the proof of \eqref{eq:3}. Concerning \eqref{eq:4}, one has on $\partial \Sigma$,

$$\begin{array}{rcl}
D_\nu^\perp Y^\perp&=&(v,D_\nu x)x^\perp+(v,\nu)D_\nu^\perp x^\perp -\frac{1}{2}(D_\nu|x|^2)v^\perp+\frac{1}{2}(1-|x|^2)D_\nu^\perp v^\perp\\\\
&=& (v,\nu)D_\nu^\perp x^\perp -\frac{1}{2}(D_\nu|x|^2)v^\perp,
\end{array}$$
since $x^\perp\equiv 0$ on $\partial \Sigma$ by the free boundary condition and $|x|\equiv1$ on $\partial\Sigma$. Using (ii) and (iv) of Lemma \ref{lem:tech}, one gets, 

$$\begin{array}{rcl}
D_\nu^\perp Y^\perp &=&-(v,\nu)\sigma(\nu,\nu)-v^\perp\\\\
&=&D_\nu^\perp v^\perp-v^\perp,
\end{array}$$
which proves \eqref{eq:4}.

\end{proof}

\noindent {\em Proof of Lemma \ref{lem:tech}:}

\medskip

The assertion (iii) follows immediately from the fact that $\Delta_\Sigma x=0$ (since $\Sigma$ is minimal). For (i), one computes

$$\begin{array}{rcl}
D_ix^\perp&=&D_i\left(\sum_{j=k+1}^n(x,N_j)N_j\right)\\\\
&=&\sum_{j=k+1}^n \left\{(D_ix,N_j)N_j+(x,D_iN_j)N_j+(x,N_j)D_iN_j\right\}\\\\
&=&\sum_{j=k+1}^n (e_i,N_j)N_j-\sigma(x^\top,e_i)+\sum_{j=k+1}^n(x,N_j)D_iN_j\\\\
&=&-\sigma(x^\top,e_i)+\sum_{j=k+1}^n(x,N_j)D_iN_j.
\end{array}$$
Notice that by assumption, $D_i^\perp N_j=0$, hence by projecting orthogonally onto the normal bundle,

$$D_i^\perp x^\perp=-\sigma(x^\top,e_i),$$
which proves (i). For (ii), a similar computation leads to

$$\begin{array}{rcl}
D_iv^\perp &=& -\sigma(v^\top,e_i)+\sum_{j=k+1}^n(v,N_j)D_iNj,
\end{array}$$
which implies as above that

$$D_i^\perp v^\perp=-\sigma(v^\top,e_i),$$
and (ii) is proved. Equation (iv) is an immediate consequence of (i) and the fact that $x=\nu$ is tangent to $\Sigma$ at the boundary. For (v), choose the basis $(e_1,\cdots,e_k)$ so that $e_k=\nu$. Then, by (ii),

$$D_\nu^\perp v^\perp=-\sigma(v^\top,\nu).$$
On the other hand, for $i<k$,

$$\begin{array}{rcl}
\sigma(e_i,\nu)&=&\sum_{j=k+1}^n(D_{e_i}x,N_j)N_j\\\\
&=&\sum_{j=k+1}^n(e_i,N_j)N_j\\\\
&=&0.
\end{array}$$
This implies that

$$\sigma(v^\top,\nu)=(v,\nu)\sigma(\nu,\nu),$$
therefore

$$D_\nu^\perp v^\perp=-(v,\nu)\sigma(\nu,\nu),$$
proving (v).

\cqfd

We now present a new spectral characterization of free boundary minimal surfaces in $\B^3$ with index $4$, based on the spectrum of the non-local operator $\mathscr{A}$. By \cite[Prop. 8.1]{FS2} (see also \cite[Cor. 7.2]{D}), if $\Sigma$ is a free boundary minimal surface in $\B^3$ with index $4$, then $\lambda_0$, the first eigenvalue of the Jacobi operator with Dirichlet boundary conditions, has to be zero. Furthermore, a first eigenfunction that is positive in the interior of $\Sigma$ is 

$$\xi=(x,N).$$
Thus, $J\xi=0$ and $\xi|_{\partial\Sigma}\equiv0$. As a interesting side remark, we note that since $\xi$ is $Q$-orthogonal to every $J$-harmonic function (as an easy consequence of Green's theorem), it follows immediately from Theorem \ref{pseudo} that

$$\int_\Sigma \xi N=0.$$
We now recall a result concerning the solution of the Dirichlet problem for the Jacobi operator, whose proof follows along the lines of the proof of \cite[Lemma 4.1]{D}, and thus will be omitted:

\begin{Lem}\label{Dirichlet}

Assume that $\lambda_0=0$ and let $\xi$ be a first eigenfunction of the Jacobi operator with Dirichlet boundary conditions. Let $u\in C^\infty(\partial\Sigma)$. Then, the Dirichlet boundary value problem:

$$\left\{\begin{array}{lcl}
J\hat{u}&=&0\mbox{ on }\Sigma\\
\hat{u}|_{\partial\Sigma}&=&u\mbox{ on }\partial\Sigma
\end{array}\right.$$
is solvable, if and only if $\int_{\partial\Sigma}u\frac{\partial \xi}{\partial\nu}=0$.

\end{Lem}
We are now prepared to state and prove our second main result, in which minimal surfaces with free boundary in the unit $3$-ball are characterized by their Dirichlet spectrum $(\lambda_k)_{k\in \mathbb{N}}$ of the Jacobi operator and spectrum $(\mu_k)_{k\in \mathbb{N}}$ of the operator $\mathscr{A}$ introduced in Section 5:

\begin{Thm}\label{charac4}

Let $\Sigma\hookrightarrow \mathbb{B}^3$ be an orientable, free boundary minimal surface that is not a disk. Then, $\Sigma$ has index $4$ if and only if $\lambda_0=0$, $\mu_0=\mu_1=\mu_2=-2$ and $\mu_3\geq0$.

\end{Thm}



As a consequence of Theorem \ref{charac4}, one obtains the following lower bound for the index form, which can be seen as a (weak) free boundary analog of the inequality \eqref{ineq_Ur} for closed minimal surfaces of $S^3$ with index $5$:

\begin{Cor}\label{good-ineq}

Let  $\Sigma\hookrightarrow \mathbb{B}^3$ be an orientable, free boundary minimal surface with index $4$. Then, for every function $u\in C^\infty(\Sigma)$ such that $Ju=0$,

$$Q(u)\geq -2\int_\Sigma |u|^2.$$

\end{Cor}

\begin{Rem}

{\em

One should note that {\em no boundary condition} whatsoever is imposed on $u$ in Corollary \ref{good-ineq}. Also, from the point of view of Calculus of Variations, the condition $Ju=0$ is quite natural: it is precisely saying that $uN$ is a Jacobi vector field on $\Sigma$, which implies that it generates a $1$-parameter family of deformations of $\Sigma$ that, it one neglects the boundary terms and the free boundary condition, preserves the area up to order two near the surface $\Sigma$.

}
\end{Rem}
Note that if one denotes (as in Lemma \ref{Dirichlet}) by $\hat{u}$ the unique $J$-harmonic extension of a function $u\in C^\infty(\partial\Sigma)$, then the infimum, over all functions $u\in C^\infty(\partial \Sigma)$ such that $\int_{\partial \Sigma}u\frac{\partial \xi}{\partial\nu}=0$, of the Rayleigh quotient

$$\frac{Q(\hat{u})}{\int_{\partial \Sigma} u^2},$$
is equal to $\sigma_1-1$, where $\sigma_1$ first eigenvalue of the Jacobi operator for the Steklov boundary value problem. As we have mentioned before, the numerical value of $\sigma_1$ does not appear to have any significant geometric meaning. On the contrary, the infimum of the alternative Rayleigh quotient

$$ \frac{Q(\hat{u})}{\int_{\Sigma} \hat{u}^2},$$
is equal to $-2$, provided the index of $\Sigma$ is equal to $4$. Furthermore, in this case the infimum is achieved with multiplicity $3$ (each component of the unit normal being a minimizer).

\begin{proof}

In order to prove each implication in Theorem \ref{charac4}, we will need the following lemma, whose proof is postponed to the end of the proof of the theorem:

\begin{Lem} \label{boundaryintegral}

The function $\xi$ satisfies:

$$Q(\mathbf{1},\xi)\neq 0.$$

\end{Lem}
Now, assume first that $\Sigma$ is a free boundary minimal surface in $\B^3$ that is not a disk. Since $\Sigma$ is not a disk, the linear subspace of $C^\infty(\Sigma)$ of functions $\{v^\perp\,;\,v\in \R^3\}$ has dimension $3$. Thus, by Theorem \ref{pseudo}, $-2$ is eigenvalue of $\mathscr{A}$ of multiplicity at least $3$. Assume by contradiction that there exists $\psi\in C^\infty(\Sigma)$ such that $\psi$ is eigenfunction of $\mathscr{A}$ associated with an eigenvalue $\mu<0$, and such that

$$Q(\psi,v^\perp)=0,\quad,\forall v\in \mathbb{R}^3.$$
For $v\in\R^3$, define

$$\tilde{v}^\perp=v^\perp-\frac{Q(\mathbf{1},v^\perp)}{Q(\mathbf{1},\xi)}\cdot \xi.$$
and similarly

$$\tilde{\psi}=\psi-\frac{Q(\mathbf{1},\psi)}{Q(\mathbf{1},\xi)}\cdot \xi.$$
Clearly,

$$Q(\tilde{v}^\perp,\mathbf{1})=Q(\tilde{\psi},\mathbf{1})=0.$$
By the Green formula, for every $J$-harmonic function $h$,

$$Q(h,\xi)=\int_{\partial\Sigma}\frac{\partial h}{\partial \nu}\xi=0.$$
In particular

$$0=Q(\xi,\xi)=Q(v^\perp,\xi)=Q(\psi,\xi),$$
which implies that

$$0=Q(\xi,\tilde{v}^\perp)=Q(\xi,\tilde{\psi}).$$
Also, 

$$Q(\tilde{\psi},v^\perp)=0,$$
and

$$Q(\tilde{v}^\perp)=Q(v^\perp)<0,\,Q(\tilde{\psi})=Q(\psi)<0,\,Q(\mathbf{1})<0.$$
Therefore, the vector space generated by $\mathbf{1}$, $\tilde{\psi}$ and $\tilde{v}^\perp$, $v\in \R^3$, is a 5-dimensional vector space on which $Q$ is negative. This contradicts the hypothesis that the index of $\Sigma$ is $4$. Thus, such a $\psi$ cannot exist, and we have proved that $\mu_0=\mu_1=\mu_2=-2$, $\mu_3\geq0$.\\

Now we prove the converse. Notice that $\mu_3\geq0$ implies that any $J$-harmonic function $h$ such that $Q(h,v^\perp)=0$ for all $v\in\R^3$ satisfies $Q(h)\geq0$. For $v\in\R^3$, define

$$\tilde{v}^\perp=v^\perp-\frac{Q(\mathbf{1},v^\perp)}{Q(\mathbf{1},\xi)}\cdot \xi,$$
so that $Q(\tilde{v}^\perp,\mathbf{1})=0$ and $Q(\tilde{v}^\perp)=Q(v^\perp)<0$. Hence, $Q$ is negative definite on the space $\mathcal{W}$, defined as the linear span of $\mathbf{1}$ and $\tilde{v}^\perp$, $v\in\R^3$. Note that $\mathcal{W}$ has dimension $4$ if $\Sigma$ is not a flat disk.

If $u\in C^\infty(\Sigma)$ is such that $Q(u,\xi)=\int_{\partial\Sigma}u\frac{\partial \xi}{\partial\nu}=0$, then by Lemma \ref{Dirichlet}, one can write $u=g+h$ with $Jh=0$ and $g|_{\partial \Sigma}=0$. By the Green formula, 

$$Q(h,g)=\int_{\partial \Sigma}g\left(\frac{\partial h}{\partial\nu}-h\right)=0,$$
so that

$$Q(u)=Q(g)+Q(h).$$
Since $Q(g,v^\perp)=Q(g,\xi)=0$ for all $v\in \R^3$, and since $Q(h,\xi)=Q(u,\xi)=0$, one has

$$Q(u,\tilde{v}^\perp)=Q(h,v^\perp).$$
By assumption $\lambda_0\geq0$, thus $Q(g)\geq0$, and the fact that $\mu_3\geq0$ then implies that $Q(u)\geq0$, provided $Q(u,\tilde{v}^\perp)=0$ for all $v\in\R^3$, and $Q(u,\xi)=0$.

Since $Q$ is negative definite on $\mathcal{W}$, one has the orthogonal decomposition:

$$W^{1,2}(\Sigma)=\mathcal{W}\oplus\mathcal{W}^\perp,$$
where $\mathcal{W}^\perp$ denotes the space of functions that are $Q$-orthogonal to $\mathcal{W}$. We claim that $Q$ is non-negative on $\mathcal{W}^\perp$, which implies that the index of $Q$ is $4$. Assume this is not true, then one finds $\psi$ in $\mathcal{W}^\perp$, such that $Q(\psi)<0$. Let 

$$\varphi=\psi-\frac{Q(\psi,\xi)}{Q(\mathbf{1},\xi)}\cdot \mathbf{1}.$$
Since $Q(\psi,\tilde{v}^\perp)=0$ and $Q(\mathbf{1},\tilde{v}^\perp)=0$, one finds that

$$Q(\varphi,\tilde{v}^\perp)=0,\,\,\forall v\in\R^3.$$
Therefore, one has $Q(\varphi)\geq0$. But, since $Q(\psi,\mathbf{1})=0$,

$$Q(\varphi)=Q(\psi)+\alpha^2Q(\mathbf{1})<0,$$
with $\alpha=\frac{Q(\psi,\xi)}{Q(\mathbf{1},\xi)}$, which is a contradiction.

\end{proof}
{\em Proof of Lemma \ref{boundaryintegral}}\\

Since $\xi>0$ in the interior of $\Sigma$ and $\xi$ vanishes on $\partial\Sigma$, one has 

$$\frac{\partial\xi}{\partial\nu}\leq 0.$$
It is thus enough to show that $\frac{\partial\xi}{\partial\nu}< 0$ in at least one point of the boundary. One computes:

$$\begin{array}{rcl}
\frac{\partial \xi}{\partial\nu} &=& D_x(x,N)\\\\
&=& (D_xx,N)+(x,D_xN).
\end{array}$$
But $Dx=id$, so $D_xx=x$, and $(D_xx,N)=0$ on the boundary. Thus,

$$\frac{\partial \xi}{\partial\nu}=(x,D_xN)=-h_{xx},$$
the coefficient of the second fundamental form of $\Sigma$ in the direction $(x,x)$. Let $t$ be a unit tangent vector to $\partial\Sigma$, then the coefficient $(t,x)$ of the second fundamental form of $\Sigma$ at the boundary is

$$h_{tx}=(D_tx,N)=(t,N)=0.$$
By minimality, it follows that on $\partial\Sigma$, the matrix of $DN$ in the basis $(x,t)$ is

$$\left(\begin{array}{rcl}
-h_{xx} & 0\\
0 & h_{xx}
\end{array}\right)$$
By contradiction, let us assume that $h_{xx}\equiv0$ on $\partial\Sigma$. It follows that $DN\equiv 0$ on $\partial\Sigma$, and thus $N$ is constant on every connected component of $\partial \Sigma$. By minimality of $\Sigma$, $N: \Sigma\to S^2$ is anti-holomorphic, so $N$ must be constant everywhere on $\Sigma$ and $\Sigma$ is planar, contradiction.

\cqfd

\begin{center}{\bf Acknowledgments} \end{center}
The author is thankful to A. Fraser for her interest and support, as well as for many stimulating discussions.


\bibliographystyle{alpha}

\end{document}